\documentclass[12pt]{amsart}
\usepackage{latexsym}
\usepackage{epsfig}
\title{Contractible open 3-manifolds with free covering translation groups}
\author{Robert Myers}
\address{Department of Mathematics, Oklahoma State University,   
Stillwater, OK 74078}
\email{myersr@@math.okstate.edu}
\thanks{Reseach at MSRI is supported in part by NSF grant DMS-9022140.}

\newtheorem{prop}{Proposition}[section]
\newtheorem{thm}[prop]{Theorem}
\begin{document}

\begin{abstract} This paper concerns the class of contractible open 
3-manifolds which are ``locally finite strong end sums'' of eventually 
end-irreducible 
Whitehead manifolds. It is shown that whenever a 3-manifold in this class 
is a covering space of another 3-manifold the group of covering translations 
must be a free group. It follows that such a 3-manifold cannot cover a 
closed 3-manifold. For each countable free group a specific uncountable 
family of irreducible open 3-manifolds is constructed whose fundamental 
groups are isomorphic to the given group and whose universal covering 
spaces are in this class and are pairwise non-homeomorphic. \end{abstract}

\maketitle
 
\section{Introduction}

Suppose $M$ is a closed, connected, orientable, irreducible 
3-manifold such that $\pi_1(M)$ is infinite. The ``universal 
covering conjecture'' states that the universal covering 
space $\widetilde{M}$ of $M$ must be homeomorphic to 
$\mathbf{R}^3$. It is known that $\widetilde{M}$ is an 
irreducible, contractible, open 3-manifold \cite{MSY}. A 
\textbf{Whitehead manifold} is an irreducible, contractible, 
open 3-manifold which is not homeomorphic to $\mathbf{R}^3$. 
The universal covering conjecture is equivalent to the 
statement that Whitehead manifolds cannot cover closed 
3-manifolds. In \cite{My genus one} the author proved that ``genus one'' 
Whitehead manifolds cannot non-trivially cover other 
3-manifolds, even non-compact ones. Wright \cite{Wr} extended this 
result to the much larger class of ``eventually 
end-irreducible'' Whitehead manifolds, a class which includes all those 
Whitehead manifolds which are monotone unions of cubes with a bounded 
number of handles. Tinsley and Wright \cite{Ti-Wr} 
gave specific examples of Whitehead manifolds which are not 
eventually end-irreducible and cannot non-trivially cover any 
other 3-manifolds. They also constructed an uncountable family 
of Whitehead manifolds which are infinite cyclic covering 
spaces of other 3-manifolds and deduced from the countability 
of the set of homeomorphism types of closed 3-manifolds that 
there must be uncountably many of these which cannot cover 
closed 3-manifolds; however their methods did not establish 
which ones these were. In \cite{My cover} the author constructed a different 
uncountable family of Whitehead manifolds which are infinite 
cyclic covering spaces of other 3-manifolds and used different 
techniques to prove that none of them covers a closed 
3-manifold. 

This paper combines the methods of \cite{My cover}, 
\cite{Wr}, and \cite{Ti-Wr} to give 
a much larger class than in \cite{My cover} of specific Whitehead manifolds 
which do not cover closed 3-manifolds but may non-trivially 
cover other non-compact 3-manifolds, namely the class of 
``strong end sums along a locally finite tree''  of eventually 
end-irreducible Whitehead manifolds. In fact it is shown that 
whenever such a manifold covers a 3-manifold the group of 
covering translations must be a free group (Theorem 3.1). 
Moreover for any countable free group there are uncountably 
many specific examples of orientable, irreducible open 
3-manifolds whose fundamental groups are isomorphic to the 
given group and whose universal covering spaces belong to 
this class and are pairwise non-homeomorphic (Theorem 4.1). 
There are also uncountably many specific examples in this 
class which can be only infinite cyclic covering spaces of 
3-manifolds and uncountably many specific examples which 
cannot non-trivially cover any 3-manifold. 

The results of \cite{My cover} use a theorem of Geoghegan and Mihalik 
\cite{Ge-Mi} 
which implies that whenever a Whitehead manifold $W$ covers 
an orientable 3-manifold the group of covering translations 
must inject into the mapping class group of $W$. If $W$ covers 
a closed, orientable, irreducible 3-manifold then the group of 
covering translations must be finitely generated and 
torsion-free. 
In \cite{My cover} the examples were constructed so that every finitely 
generated, torsion-free subgroup of their mapping class 
groups must have a subgroup of finite index which either has 
infinite abelianization or a non-trivial normal abelian 
subgroup. Results of Waldhausen \cite{Wa}, Hass-Rubinstein-Scott 
\cite{HRS}, Mess \cite{Me}, 
Casson-Jungreis \cite{Ca-Ju}, and Gabai \cite{Ga} were then quoted which 
imply that a closed, orientable, irreducible 3-manifold with such a 
fundamental group must have universal covering space 
homeomorphic to $\mathbf{R}^3$. 

The present paper avoids the use of the Geoghegan-Mihalik 
result and the requisite analysis of the mapping class group. 
For the class of Whitehead manifolds under consideration 
results of \cite{My endsum} are used to show that the group of covering 
translations acts on a certain simplicial tree. The Orbit 
Lemma of \cite{Wr} and the Special Ratchet Lemma of \cite{Ti-Wr} are then 
used to prove that this action fixes no vertices, from which the 
result follows. We remark that the methods by which Tinsley and Wright 
apply these lemmas in their proof of Theorem 5.3 of \cite{Ti-Wr} 
could be adapted to prove this fact. However, we present a different, 
somewhat more direct argument which is closer in spirit to Wright's 
proof of the main theorem of \cite{Wr}. We also give an alternative, 
somewhat shorter proof of the special case of the Orbit Lemma that we 
use.

The Whitehead manifolds considered in \cite{Ti-Wr}, \cite{My cover}, and this 
paper are all ``end sums'' of Whitehead manifolds; they are obtained by 
gluing together a collection of Whitehead manifolds in a certain way (see the 
next section for the precise definition). The summands in \cite{Ti-Wr} are 
members of a certain uncountable collection of genus one 
Whitehead manifolds discovered by McMillan \cite{Mc}; the summands in 
\cite{My cover} are members of a different uncountable collection of genus 
one Whitehead manifolds chosen so that the mapping class group of the end sum 
will have the appropriate structure as described 
above. However, the main difference is not in the summands, but in how they 
are glued together. The examples of \cite{My cover} are all ``strong'' end 
sums which have a certain ``rigidity up to isotopy'' in their construction. 
The end sums in \cite{Ti-Wr} are not strong end sums; in fact it follows 
from Proposition 2.1 below that these  manifolds cannot be expressed in any 
way as strong end sums, even though by Proposition 2.2 below their summands 
can be glued together in a different fashion to obtain different manifolds 
which are strong end sums. Thus the results of this paper 
apply to all the examples of \cite{My cover} but none of the examples of 
\cite{Ti-Wr}. The question of which of them cannot cover closed 
3-manifolds (conjecturally all of them) is still open. 

\section{Background Material} 

For general background on 3-manifold topology see \cite{He} or \cite{Ja}. 
We denote the manifold theoretic boundary and interior of a 
manifold $M$ by $\partial M$ and $int \, M$, respectively. 
We denote the topological boundary, interior, and closure of 
a submanifold $M$ of a manifold $N$ by $Fr_N(M)$, $Int_N(M)$, 
and $Cl_N(M)$, respectively, with the subscript deleted when 
$N$ is clear from the context. The \textbf{exterior} of 
$M$ in $N$ is the closure of the complement of a regular 
neighborhood of $M$ in $N$. $M$ is \textbf{open} if 
$\partial M=\emptyset$ and no component of $M$ is compact. 
A continuous map $f:M \rightarrow N$ of manifolds is 
\textbf{$\partial$-proper} if $f^{-1}(\partial N)=\partial M$. 
It is \textbf{end-proper} if preimages of compact sets are 
compact. It is \textbf{proper} if it has both these properties. 
These terms are applied to a submanifold if its inclusion map 
has the corresponding property. Two codimension one 
submanifolds $M_0$ and $M_1$ of $N$, each of which is either 
proper in $N$ or is a submanifold of $\partial N$, are 
\textbf{parallel} if some component of $N-(M_0 \cup M_1)$ 
has closure homeomorphic to $M_0 \times [0,1]$ with 
$M_i=M_0 \times \{i\}$, $i=0,1$. A proper codimension one 
submanifold of $N$ is \textbf{$\partial$-parallel} if it is 
parallel to a submanifold of $\partial N$. 

An \textbf{exhaustion} $\{K_n\}_{n\geq 0}$ for a connected, 
non-compact manifold $W$ is a sequence of compact, connected, 
codimension zero submanifolds of $W$ whose union is $W$, 
such that $K_n \subseteq Int \, K_{n+1}$, $K_n \cap \partial 
W$ is either empty or a codimension zero submanifold of 
$\partial W$, and $W-Int \, K_n$ has no compact components. 

A connected, non-compact 3-manifold $W$ is \textbf{eventually 
end-irreducible} if it has an exhaustion $\{K_n\}$ such that 
$Fr \, K_n$ is incompressible in $W-Int \, K_0$ for all 
$n \geq 0$. We also say that $W$ is \textbf{end-irreducible 
rel $K_0$}. $W$ is \textbf{eventually $\pi_1$-injective at 
$\infty$} if there is a compact subset $J$ of $W$ such that 
for every compact subset $K$ of $W$ containing $J$ there is 
a compact subset $L$ of $W$ containing $K$ such that every 
loop in $W-L$ which is null-homotopic in $W-J$ is 
null-homotopic in $W-K$. We also say that $W$ is 
\textbf{$\pi_1$-injective at $\infty$ rel $J$}. It is a 
standard exercise to show that $W$ is eventually 
end-irreducible if and only if it is eventually 
$\pi_1$-injective at $\infty$. Note in particular that if $W$ 
is end-irreducible rel $K_0$, then it is $\pi_1$-injective at 
$\infty$ rel $K_0$. 

Let $V$ be an irreducible non-compact 3-manifold such that either 
$\partial V=\emptyset$ 
or each component of $\partial V$ is a plane. A proper plane $P$ in $V$ is 
\textbf{trivial} if some component of $V-P$ has closure homeomorphic to 
$\mathbf{R}^2 \times [0,\infty)$ with $P=\mathbf{R}^2 \times \{0\}$. 
$V$ is \textbf{$\mathbf{R}^2$-irreducible} every 
proper plane in $V$ is trivial (hence $\partial V=\emptyset$ or 
$V=\mathbf{R}^2 \times [0,\infty))$; it is \textbf{aplanar} if every  
proper plane in $V$ is either trivial or $\partial$-parallel. 
A \textbf{partial plane} is a simply connected, non-compact 2-manifold 
with non-empty boundary. $V$ is \textbf{strongly aplanar} if it is 
aplanar and given any proper 2-manifold $\mathcal{P}$ in $V$ each 
component of which is a partial plane, there is a collar on $\partial V$ 
which contains $\mathcal{P}$. $V$ is \textbf{anannular at $\infty$} if 
for every compact subset $K$ of $V$ there is a compact subset $L$ of $V$ 
containing $K$ such that $V-L$ is \textbf{anannular}, i.e. every proper, 
incompressible annulus in $V-L$ is $\partial$-parallel. 

Now suppose we are given a countable simplicial tree $\Gamma$ to each 
vertex $v_i$ of which we have associated a connected, oriented, irreducible, 
non-compact 3-manifold $V_i$ whose boundary is a non-empty disjoint 
union of planes. Suppose that to each edge $e_k$ of $\Gamma$ we have 
associated a component of $\partial V_i$ and a component of $\partial V_j$, 
where $e_k$ has endpoints $v_i$ and $v_j$ and no boundary plane is associated 
to different edges. The connected, oriented, non-compact 3-manifold $W$ 
obtained by gluing each such pair of planes by an orientation reversing 
homeomorphism is called the \textbf{plane sum} of the $V_i$ along $\Gamma$. 
The image in $W$ of the pair of planes identified as above is denoted by 
$E_k$ and is called a \textbf{summing plane}. The plane sum is 
\textbf{degenerate} if either some summing plane is trivial or 
$\partial$-parallel in $W$ or two distinct summing planes are parallel 
in $W$. Theorem 3.2 of \cite{My endsum} gives necessary and sufficient 
conditions on 
$\Gamma$ and the $V_i$ for the plane sum to be non-degenerate.  For our 
present purposes Corollary 3.3 of \cite{My endsum}, which states that 
the plane sum is 
non-degenerate if no summand $V_i$ has a boundary plane $E_k$ such that 
$E_k \cup int \, V_i$ is homeomorphic to $\mathbf{R}^2 \times [0,\infty)$, 
will suffice because in our case $int \, V_i$ will be a Whitehead manifold.  
The plane sum is \textbf{strong} if it is non-degenerate and 
each summand is strongly aplanar and anannular at $\infty$. 

\begin{prop} Let $W$ be a non-degenerate plane sum of aplanar 3-manifolds 
along a locally finite tree. Let $W^{\prime}$ be a strong plane sum. 
Let $\mathcal{E}$ and $\mathcal{E}^{\prime}$ be the unions of the 
respective sets of summing planes. Suppose $g:W \rightarrow W^{\prime}$ 
is a homeomorphism. Then $g$ is ambient isotopic rel $\partial W$ to 
a homeomorphism $h$ such that $h(\mathcal{E})=\mathcal{E}^\prime$. 
\end{prop}

\noindent \textbf{Proof:} This is Theorem 4.3 of \cite{My endsum}. $\Box$

\bigskip

Now suppose that given $\Gamma$ we have associated to each vertex $v_i$ 
a connected, open, irreducible, oriented 3-manifold $W_i$, and that to each 
edge $e_k$ we have associated an end-proper ray (a space homeomorphic to 
$[0,\infty)$) in $W_i$ and an end-proper 
ray in $W_j$, where $e_k$ has endpoints $v_i$ and $v_j$, the rays associated 
to different edges are disjoint and their union is end-proper. The exterior 
$V_i$ of the union of the rays contained in $W_i$ is then bounded by planes. 
Note that $int \, V_i$ is homeomorphic to $W_i$. The plane sum $W$ of the 
$V_i$ along $\Gamma$ is called an \textbf{end sum} of the $W_i$ along 
$\Gamma$. (Note that $W$ depends on the choice of the rays; this dependence 
is investigated further in \cite{My endsum}.) A \textbf{strong end sum} is 
one whose 
associated plane sum is strong. 

We conclude this section with some remarks about the existence of strong 
end sums. In the present context the following is the most relevant fact; 
more general results may be found in \cite{My attach} and \cite{My endsum}. 

\begin{prop} Given a countable, locally finite tree $\Gamma$, a  
collection $\{W_i\}$ of connected, irreducible, oriented, one-ended open 
3-manifolds, and a bijection between the vertices of $\Gamma$ and $\{W_i\}$, 
there exists a strong end sum of the $W_i$ along $\Gamma$. \end{prop}

\noindent \textbf{Proof:} This is a special case of Theorem 5.1 of 
\cite{My endsum}. $\Box$

\bigskip

For later reference we briefly describe the construction of the rays 
required in the proof of this result. Suppose $V$ is a connected, orientable, 
irreducible, one-ended, non-compact 3-manifold whose boundary is either 
empty or consists of a finite set of disjoint planes. An exhaustion 
$\{C_n\}$ for $V$ is \textbf{nice} if for all $n \geq 1$ one has that 
$C_n-Int \, C_{n-1}$ is irreducible, $\partial$-irreducible, and anannular, and that 
for all $n \geq 0$ one has that 
each component of $Fr \, C_n$ has positive genus and negative Euler 
characteristic, and if $\partial V 
\neq \emptyset$, one has that $C_n \cap \partial V$ consists of a single 
disk in each component of $\partial V$. One says that $V$ is \textbf{nice} 
if it has a nice exhaustion. 

\begin{prop} If $V$ is nice, then $V$ is strongly aplanar and anannular 
at $\infty$. \end{prop}

\noindent \textbf{Proof:} This follows from Theorem 5.3 and Lemma 1.3 (6) 
of \cite{My attach}. $\Box$

\bigskip

Given $W_i$ one chooses an exhaustion $\{K_n\}$ for $W_i$ with each 
$\partial K_n$ connected and of positive genus. If $\nu$ rays are required, 
then for each $n\geq 1$  
one chooses a disjoint union of $\nu$ proper arcs in $K_n-Int \, K_{n-1}$ 
each component of which joins $Fr \, K_{n-1}$ to $Fr \, K_n$. This is 
done so that the endpoints match up on $Fr \, K_n$ so as to give 
$\nu$ rays in $W_i$. Then we obtain an exhaustion $\{C_n\}$ for $V_i$ by 
letting $C_0=K_0$ and for $n\geq 1$ letting $C_n$ be the exterior in 
$K_n-Int \, K_{n-1}$ of the union of the arcs. All that remains is to note 
that by Theorem 1.1 of \cite{My excel} one can choose the arcs so that 
$C_n-Int \, C_{n-1}$ is irreducible, $\partial$-irreducible, and anannular. 

In section 4 we will give explicit constructions of examples of this type 
which do not rely on Theorem 1.1 of \cite{My excel}. 

\section{The General Result}

\begin{thm} Let $W$ be a strong end sum of eventually end-irreducible 
Whitehead manifolds $W_i$ along a locally finite tree $\Gamma$. If $W$ is a 
covering space of a 3-manifold $M$, then there is a simplicial action of 
$\pi_1(M)$ on $\Gamma$ under which no non-trivial element of $\pi_1(M)$ 
fixes a vertex of $\Gamma$. Hence 
\begin{enumerate}

\item $\pi_1(M)$ is a free group. 

\item $M$ cannot be a closed 3-manifold. 

\item If $\Gamma$ has countably many ends, then $\pi_1(M)$ is cyclic. 

\item If the number of ends of $\Gamma$ is finite and greater than 
two, then $\pi_1(M)$ is trivial, i.e. $M=W$.  
\end{enumerate} 
\end{thm}

\noindent \textbf{Proof:} We first show how to deduce (1)--(4) from the 
main statement of the theorem. (1) $\pi_1(M)$ has a subgroup $H$ of index at 
most two which acts on $\Gamma$ without inversions of the edges, hence 
acts freely on $\Gamma$, hence is free. It follows that $\pi_1(M)$ is 
itself free \cite{St}. (2) If $M$ were closed then it would be a connected 
sum of 
2-sphere bundles over $S^1$ \cite[Theorem 5.2]{He}, hence would not be 
aspherical, hence its 
universal covering space would not be contractible. (3) If $rank \, 
\pi_1(M) \geq 2$, then $\Gamma$ has uncountably many ends. (4) Suppose $A$ is 
an axis for the action of $\pi_1(M)$ on $\Gamma$, i.e. $A$ is a subtree 
isomorphic to a triangulation of $\mathbf{R}$ which is invariant under 
the infinite cyclic action (see \cite{Sh}). Since $\Gamma$ has at least three 
ends some component of $\Gamma-A$ has non-compact closure $T$, and the 
translates of $T$ yield infinitely many ends of $\Gamma$. 

We now prove the main statement of the theorem. Let $G \cong \pi_1(M)$ 
be the group of covering translations. By Proposition 2.1 each $g \in G$ 
is isotopic to a homeomorphism $h$ such that 
$h(\mathcal{E})=\mathcal{E}$, where $\mathcal{E}$ is the union 
of the summing planes of $W$. Thus $h$ determines an element of 
$Aut(\Gamma)$. We claim that this element depends only on $g$. We repeat 
the argument of Theorem 3.2 of \cite{My cover}. If $h^{\prime}$ were a homeomorphism 
isotopic to $g$ which determined a different automorphism, then $h$ and 
$h^{\prime}$ would send some summing plane $E_i$ to different summing 
planes $E_j$ and $E_k$, hence they would be ambient isotopic in $W$. 
But by Theorem 5 of \cite{Wi} disjoint, ambient isotopic, non-trivial, proper 
planes in an irreducible 3-manifold must be parallel. This contradicts 
the non-degeneracy of strong end sums. 

Thus we have a well defined action of $G$ on $\Gamma$. We next state the 
results of \cite{Wr} and \cite{Ti-Wr} that we shall need in order to prove that no vertex 
is fixed by a non-trivial element of $G$. 

Let $G$ be a group acting on an $n$-manifold $W$. One says that $G$ acts 
\textbf{without fixed points} if the only element of $G$ fixing a 
point is the identity. $G$ acts \textbf{totally discontinuously} 
if for every compact 
subset $C$ of $W$ one has that $g(C) \cap C = \emptyset$ for all 
but finitely many elements of $G$. (In \cite{Wr} the term  
``properly discontinuously" is used for this property; 
we follow Freedman and Skora's terminology \cite{Fr-Sk} in order to avoid 
confusion with other meanings of this term.) Let $p:W \rightarrow Y$ be 
the projection to the orbit space $Y$ of the action. Then $G$ acts without 
fixed points and totally discontinuously on $W$ if and only if $p$ is a 
regular covering map with group of covering translations $G$ and $Y$ is 
an $n$-manifold. (See \cite{Ma}.) In this case if $W$ is contractible, 
then $G$ must be torsion-free (see e.g. \cite{My genus one} or \cite{Wr}). 

\begin{prop}[Orbit Lemma (Wright)] Let $W$ be a con\-tract\-i\-ble, open 
$n$-man\-i\-fold, $n \geq 3$. Let $g$ be a non-trivial homeomorphism of $W$ 
onto itself such that the group $<g>$ of homeomorphisms generated by 
$g$ acts without fixed points and totally discontinuously on $W$. 
Given compact subsets $B$ and $Q$ of $W$, there is a compact subset 
$C$ of $W$ containing $B$ such that every loop in $W-C$ is homotopic in $W-B$ 
to a loop in $W-\cup_{i=-\infty}^{\infty} g^i(Q)$. \end{prop}

\noindent \textbf{Proof:} Except for the statement that $C$ contains 
$B$ this is Lemma 4.1 of \cite{Wr}; we can clearly enlarge the $C$ of that 
result to satisfy this requirement. 

We now give an alternate proof for the special case in which $W$ is an 
irreducible 3-manifold. The quotient manifold $Y=W/<g>$ is an irreducible 
open 3-manifold having the homotopy type of a circle. Any irreducible open 
3-manifold with locally free fundamental group has an exhaustion by cubes 
with handles (Theorem 2 of \cite{F-F}). Let $\{Y_n\}$ be such an 
exhaustion for $Y$. We may assume that $\pi_1(Y_0) \rightarrow \pi_1(Y)$ 
is onto and $p(Q) \subseteq Int \, Y_0$, where $p:W \rightarrow Y$ is the 
covering projection. Thus $\cup^{\infty}_{i=-\infty} g^i(Q) \subseteq 
Int \, p^{-1}(Y_0)$. Now $p^{-1}(Y_0)$ is a non-compact cube with handles. 
There is a finite set of disjoint, proper disks in $p^{-1}(Y_0)$ whose 
union splits $p^{-1}(Y_0)$ into a compact cube with handles $H$ which 
contains $B \cap p^{-1}(Y_0)$ and a 3-manifold $H^{\prime}$ whose 
components are non-compact cubes with handles. These splitting disks can 
be chosen disjoint from $B$. Let $C=B \cup H$. Suppose $\gamma$ is a loop 
in $W-C$. Homotop $\gamma$ so that it is in general position with respect 
to $\partial H^{\prime}$. Then it meets $H^{\prime}$ in a finite set of 
paths $\gamma_j$. Since  the components of $H^{\prime}$ are cubes with 
handles each $\gamma_j$ can be homotoped rel $\partial \gamma_j$ 
to a path $\gamma_j^{\prime}$ in $\partial H^{\prime}$.  
This can be done so that no 
$\gamma_j^{\prime}$ meets a splitting disk. Thus $\gamma$ is homotopic in 
$W-C$, and hence in $W-B$, to a loop $\gamma^{\prime}$ which lies in 
$W-Int \, p^{-1}(Y_0)$ and hence in $W-\cup^{\infty}_{i=-\infty} g^i(Q)$. 
$\Box$ 

\bigskip 

\begin{prop}[Special Ratchet Lemma (Tinsley-Wright)] Let $W$ be an 
op\-en $n$-manifold and $W_0$ an open subset of $W$ with closure $V_0$. 
Suppose $W_0$ is $\pi_1$-injective at $\infty$ rel $J$, $V_0$ is an 
$n$-manifold, $\partial V_0$ is proper and bicollared in $W$, and each 
component of $\partial V_0$ is simply connected. Let $g$ be a homeomorphism 
of $W$ onto itself such that each of $g(J)$ and $g^{-1}(J)$ can be 
ambiently isotoped into $W_0$. Then there is a compact subset $R$ of $W$ 
containing $J$ 
such that a loop in $W-\cup_{i=-\infty}^{\infty} g^i(R)$ is null-homotopic 
in $W-J$ if and only if it is null-homotopic in $W-g^i(J)$ for each 
$i \in \mathbf{Z}$. \end{prop}

\noindent \textbf{Proof:} This is a slight variation of Lemma 5.1 of 
\cite{Ti-Wr} which has the same proof. $\Box$ 

\bigskip

The hypotheses of the Special Ratchet Lemma are clearly satisfied when 
$G$ acts on $\Gamma$ with fixed points, i.e. some non-trivial $g \in G$ 
is isotopic to $h$ such that $h(V_0)=V_0$ for the plane summand $V_0$ 
associated to an end summand $W_0$. We shall prove that $W_0$ is 
\textbf{$\pi_1$-trivial at $\infty$}, i.e. for every compact subset 
$A$ of $W_0$ there is a compact subset $A^*$ of $W_0$ containing 
$A$ such that every loop in $W_0-A^*$ is null-homotopic in $W_0-A$. 
By a result of C. H. Edwards \cite{Ed} and 
C. T. C. Wall \cite{Wl} every irreducible, contractible, open 
3-manifold which is $\pi_1$-trivial at $\infty$ must be homeomorphic 
to $\mathbf{R}^3$. This contradicts the assumption that $W_0$ is a 
Whitehead manifold. 

So, let $A$ be a compact subset of $W_0$. Now $W_0$ is $\pi_1$-injective 
at $\infty$ rel $J$ for some compact subset $J$ of $W_0$. By the Special 
Ratchet Lemma there is a compact subset $R$ of $W$ containing $J$ such that 
a loop in $W-\cup_{i=-\infty}^{\infty} g^i(R)$ is null-homotopic in $W-J$ 
if and only if it is null-homotopic in $W-g^i(J)$ for all $i \in \mathbf{Z}$. 
Let $N=\partial V_0 \times [0,1]$ be a collar on $\partial V_0$ in $V_0$ 
such that $\partial V_0 \times \{0\}=\partial V_0$ and $N \cap (A \cup J)=
\emptyset$. Let $R_0=R \cap Cl(V_0-N)$. Then $R_0$ is a 
compact subset of $W_0$ which contains $J$. Let $K=A \cup R_0$. 
Since $W_0$ is $\pi_1$-injective at $\infty$ rel $J$ there is a compact 
subset $L$ of $W_0$ containing $K$ such that loops in $W_0-L$ which are 
null-homotopic in $W_0-J$ are null-homotopic in $W_0-K$. Apply the Orbit 
Lemma with $B=L$ and $Q=R$ to get a compact subset $C$ of $W$ containing 
$L$ such that 
every loop in $W-C$ is homotopic in $W-L$ to a loop in 
$W-\cup_{i=-\infty}^{\infty} g^i(R)$. By enlarging $C$, if necessary, we 
may assume that $C \cap N$ consists of cylinders $D_j \times [0,1]$, 
where $D_j$ is a disk in the component $E_j$ of $\partial V_0$. 
There is an $s \in (0,1)$ such that the collar $N_s=\partial V_0 \times 
[0,s]$ misses $L$. Let $C_0=C \cap Cl(V_0-N_s)$. 

We claim that we may take $A^*=C_0$. Consider a loop $\gamma$ in 
$W_0-C_0$. We will show that $\gamma$ is null-homotopic in $W_0-A$. 
First note that $\gamma \cap C$ is contained in the union of the 
$D_j \times (0,s)$. We can homotop $\gamma$ in $W_0-C_0$, if necessary, 
so that it misses the union of the $\{x_j\} \times [0,s]$, where $x_j$ 
is a point in the interior of $D_j$. By pushing radially outward from 
$\{x_j\} \times [0,s]$ in each $D_j \times [0,s]$ and then off 
$D_j \times [0,s]$ we obtain a homotopy of $\gamma$ in $W_0-C_0$ to a loop 
$\gamma^{\prime}$ in $W_0-(W_0 \cap C)$. Now $\gamma^{\prime}$ is 
homotopic in $W-L$ to a loop $\gamma^{\prime\prime}$ in 
$W-\cup_{i=-\infty}^{\infty} g^i(R)$. Since $W$ is contractible 
$\gamma^{\prime\prime}$ is null-homotopic in $W$. Since $<g>$ is totally 
discontinuous $\gamma^{\prime\prime}$ is null-homotopic in 
$W-g^i(J)$ for some $i$. Since $\gamma^{\prime\prime}$ lies in $W-L$ 
the Special Ratchet Lemma implies that $\gamma^{\prime\prime}$ is 
null-homotopic in $W-J$. Since $J \subseteq L \subseteq C_0$ we have 
that $\gamma$ is null-homotopic in $W-J$. Since $\gamma$ lies in $W_0-J$ 
and the components of $\partial V_0$ are simply connected we have that 
$\gamma$ is null-homotopic in $W_0-J$. Thus $\gamma$ is null-homotopic 
in $W_0-K \subseteq W_0-A$, as required. $\Box$ 


\section{Specific Examples}

\begin{thm} 
\hspace{1in} 

\begin{enumerate}  

\item Given any countable free group $F$ there are uncountably many 
specific irreducible, orientable, open 3-manifolds $X$ such that 
$\pi_1(X)\cong F$, any 3-manifold $M$ covered by the universal covering 
space $W$ of $X$ must have free fundamental group, and the $W$ are 
pairwise non-homeomorphic. 

\item If $F\cong \mathbf{Z}$, then $X$ can be chosen so that 
$\pi_1(M)$ must be infinite cyclic.

\item If $F$ is trivial, then $X=W$ can be chosen so that $M=W$. 
\end{enumerate} 
\end{thm}

\noindent \textbf{Proof:} (1) It suffices to consider the case when 
$F$ has rank two. The construction will be a generalization of that of 
Theorem 6.1 of \cite{My cover}. Figure 1 shows a six component tangle 
$\lambda$ in 
a 3-ball $B$ called the true lover's 6-tangle. By Proposition 4.1 of 
\cite{My homology} the exterior of 
$\lambda$ is \textbf{excellent}, i.e. it is irreducible, 
$\partial$-irreducible, anannular, and atoroidal, contains a proper 
incompressible surface, and is not a 3-ball. It follows immediately from 
the proof of the result cited that each of the $k$-tangles consisting of 
$k \geq 2$ consecutive components of $\lambda$ also has excellent exterior. 
By sliding the endpoints of the arcs of $\lambda$ one sees that the 
exterior of $\lambda$ is homeomorphic to the exterior of the graph  
$\xi$ in Figure 2. By deleting the first, second, fifth, and sixth arcs 
we obtain the 2-tangle $\mu$ in Figure 3, which thus has excellent 
exterior. 

\begin{figure}[b]
\begin{center}
\epsfig{file=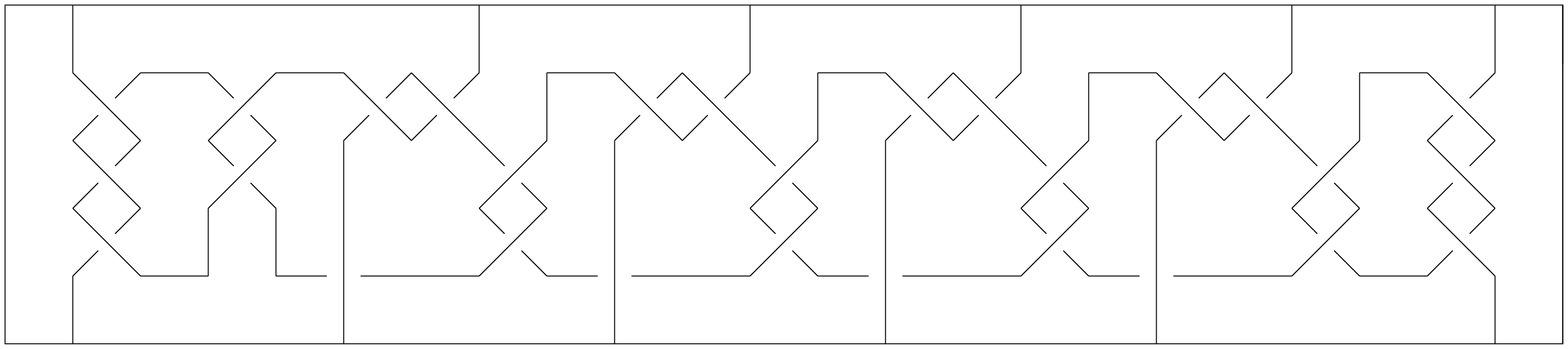, width=4in}
\end{center}
\caption{The 6-tangle $\lambda$}
\end{figure} 


\begin{figure}[b]
\begin{center}
\epsfig{file=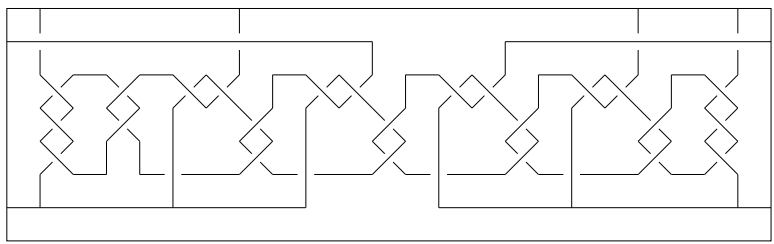, width=4in}
\end{center}
\caption{The graph $\xi$}
\end{figure}


\begin{figure}
\begin{center}
\epsfig{file=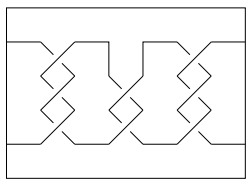, width=1.25in}
\end{center}
\caption{The 2-tangle $\mu$}
\end{figure}


We next identify the disks which are the left and right sides of the 
rectangular solid $B$ in Figures 2 and 3 
to obtain a solid torus $K$. This is done so that 
$\mu$ becomes a simple closed curve $\sigma$ and $\xi$ becomes a graph 
$\theta$ consisting of $\sigma$ together with four disjoint arcs 
$\alpha^1$, $\alpha^2$, $\alpha^3$, $\alpha^4$ joining $\sigma$ to 
$\partial K$. It follows from Lemma 2.1 of \cite{My excel} that the 
exteriors of 
$\sigma$ and of $\theta$ in $K$ are excellent. 

Now let $L$ be a regular neighborhood of $\sigma$ in $K$. We construct 
a genus one Whitehead manifold $U$ with exhaustion $\{K_n\}$ by using 
as models for $(K_n, K_{n-1})$ the pair $(K,L)$. This is done so that 
the copies $\alpha^j_n$ of the $\alpha^j$ match up along their 
endpoints to give end-proper rays $\rho_j$ in $U$. We then let $V$ be 
$U$ minus the interior of a regular neighborhood $N$ of the union of 
these rays. We choose $N$ so that its intersection $N_n$ with 
$K_n-int \, K_{n-1}$ is a regular neighborhood of the union of the 
$\alpha^j_n$. We then let $C_n$ be $Cl(K_n-N_n)$ for $n\geq 1$ and 
$C_0=K_0$. Since $K_n-int \, K_{n-1}$ and $C_n-Int \, C_{n-1}$ are 
excellent we have that $U$ is an eventually end-irreducible Whitehead 
manifold and $V$ is nice. 

We now identify the boundary planes of $V$ in pairs to obtain an 
orientable 3-manifold $X$ with $\pi_1(X)$ free of rank two. The universal 
covering space $W$ of $X$ is then an end sum of Whitehead manifolds 
$W_i$ each of which is homeomorphic to $U$ such that the plane summands 
$V_i$ are homeomorphic to $V$. We then apply Theorem 3.1. 

We next show how to get uncountably many examples of this type with 
pairwise non-homeomorphic universal covering spaces. 

If one changes the sense of the central clasp in the figures by changing 
the two overcrossings to undercrossings, thereby getting a new $\sigma$ 
and $\theta$, then the same arguments show that their exteriors in $K$ 
are excellent. Denote the old and new versions by the subscripts 0 and 
1, respectively. Embed $K$ in $S^3$ in a standard way so that a line 
segment running along the bottom front edge of $B$ becomes a simple closed 
curve $\ell$ in $\partial K$ which bounds a disk in $S^3-int \, K$. 
Then $\sigma_0$ and $\sigma_1$ become the knots $8_5$ and $8_{19}$ in $S^3$ 
with normalized Alexander polynomials 
$5-4(t+t^{-1})+3(t^2+t^{-2})-(t^3+t^{-3})$ and 
$1-(t^2+t^{-2})+(t^3+t^{-3})$, respectively. It then follows that 
there is no homeomorphism from the exterior of $\sigma_0$ in $K$ to 
that of $\sigma_1$ in $K$ which carries $\ell$ to a curve homologous to 
$\pm \ell$. since if there were, then one could extend it to a 
homeomorphism of the exteriors in $S^3$ of these two knots. 

Let $s=\{s_n\}_{n\geq 1}$ be an infinite sequence of 0's and 1's. 
Carry out the construction as before by modeling the pair $(K_n, K_{n-1})$, 
for $n \geq 1$, on $(K,L_i)$, where $L_i$ is a regular neighborhood of 
$\sigma_i$ in $K$ and $i=s_n$. Do this so that the copy $\ell_n$ of $\ell$ 
in $\partial K_n$ is null-homologous in $K_{n+1}-int \, K_n$. (Up to 
orientation and isotopy there is a unique such curve.) 

Label the various manifolds arising in the construction associated to $s$ 
by a superscript $s$. If $f:U^s \rightarrow U^t$ is a homeomorphism, 
then Lemma 3.3 of \cite{My genus one} implies that $f$ can be isotoped 
so that for some 
$a$ and $b$ one has $f(K^s_{a+m})=K^t_{b+m}$ for all $m \geq 0$. Thus 
$s_{a+m}=t_{b+m}$ for all $m \geq 0$. 

One could now note that this last equation generates an equivalence 
relation on the set $\{0,1\}^{\omega}$ of all such sequences and that 
there are uncountably many equivalence classes. In keeping with the 
desire to make our examples as explicit as possible, however, we prefer 
a more concrete approach which exhibits an explicit subset $\mathcal{S}$ 
of $\{0,1\}^{\omega}$ for which the corresponding Whitehead manifolds are 
non-homeomorphic. We define $\mathcal{S}$ and define a bijection 
$\varphi : \{0,1\}^{\omega} \rightarrow \mathcal{S}$ as follows. 
Let $x \in \{0,1\}^{\omega}$. Then $s=\varphi(x)$ will consist of strings 
of consecutive 0's which are separated by single 1's. The length of the 
$n^{th}$ string of 0's is $d_n=r_1r_2 \cdots r_n$, where 
$r_j=3^{(2^{j-1})}$ if $x_j=0$ and $r_j=5^{(2^{j-1})}$ if $x_j=1$. 
Thus $d_n=3^u5^v$, where the total exponent sum $u+v=1+2+4+8+ \cdots 
+2^{n-1}=2^n-1$. 

Suppose $t=\varphi(y)$ is another sequence such that for some $a$ and $b$ 
one has $s_{a+m}=t_{b+m}$ for all $m>0$. Locate the first 1 in this common 
tail. It is followed by a string of $3^u5^v$ 0's for some unique $u$ and 
$v$. Then $u+v=2^n-1$ for a unique $n$, and so this is the $n^{th}$ 
string of 0's in both $s$ and $t$. Note that $n>1$. Suppose $d_n=
r_1r_2 \cdots r_n=q_1q_2 \cdots q_n$ where the $r_j$ and $q_j$ correspond 
to the $x_j$ and $y_j$ as above. Then $d_{n-1}=d_n/r_n$; let  
$p_{n-1}=d_n/q_n$. If $r_n=3^{(2^{n-1})}$, then since $p_{n-1}$ has 
exponent sum in 3 at most $2^{n-1}-1$ we must have $q_n=3^{(2^{n-1})}$; 
since a similar argument holds for powers of 5 we have that $r_n=q_n$. 
We inductively conclude that $r_j=q_j$, and hence $x_j=y_j$, for 
$1\leq j\leq n$. Applying this argument to all $n^{\prime}>n$ we get that 
$x=y$ and $s=t$. 

Thus we have uncountably many non-homeomorphic genus one Whitehead 
manifolds $U^s$. We construct the corresponding $V^s$, $X^s$, and $W^s$. 
The $W^s_i$ are all homeomorphic to $U^s$. It then follows from 
Proposition 2.1 that if $W^s$ and $W^t$ are homeomorphic, so are 
$U^s$ and $U^t$, hence $s=t$. 

(2) We perform the analogous construction with the first and last arcs 
deleted. See Theorem 6.1 of \cite{My cover}. 

(3) One can carry out the construction of $V$ as above with any finite 
number $\nu$ of boundary planes by using the true lover's $\nu +2$-tangle. 
Thus given any locally finite tree $\Gamma$ one can construct the 
corresponding strong end sum. One can then choose $\Gamma$ to have the 
wrong number of ends or, for variety, let $\Gamma$ be arbitrary but choose 
one $W_0$ which is not homeomorphic to any of the other $W_i$, thereby 
creating a fixed vertex for the action on $\Gamma$. $\Box$


\begin{thebibliography}{99}

\bibitem{Ca-Ju}
A. Casson and D. Jungreis, {\it Convergence groups and Seifert fibered 3-manifolds},
Invent. Math. 118 (1994), 441--456.

\bibitem{Ed}
C. H. Edwards, {\it Open 3-manifolds which are simply connected at infinity}, 
Proc. Amer. Math. Soc. 14 (1963), 391--395. 

\bibitem{F-F} 
B. Freedman and M. Freedman, {\it Kneser-Haken finiteness for bounded 
3-manifolds, locally free groups, and cyclic covers}, preprint. 

\bibitem{Fr-Sk}
M. Freedman and R. Skora, {\it Strange actions of groups on spheres}, 
J. Differential Geometry 25 (1987), 75--98.

\bibitem{Ga}
D. Gabai, {\it Convergence groups are Fuchsian groups},
Annals of Math.  136 (1992), 447--510.

\bibitem{Ge-Mi}
R. Geoghegan and M. Mihalik, {\it The fundamental group at infinity}, 
Topology 35 (1996), 655--669.

\bibitem{HRS}
J. Hass, H. Rubinstein, and P. Scott, {\it Compactifying coverings of closed
3-manifolds}, J. Differential Geometry 30 (1989), 817--832.


\bibitem{He}
J. Hempel, {\it 3-Manifolds}, Ann. of Math. Studies, No. 86, Princeton,
(1976).

\bibitem{Ja}
W. Jaco, {\it Lectures on three-manifold topology}, CBMS Regional
Conference Series in Math., No. 43, Amer. Math. Soc. (1980).

\bibitem{Ma}
W. Massey, {\it Algebraic Topology: An Introduction}, Graduate Texts in 
Mathematics No. 56, Springer-Verlag (1977). 

\bibitem{Mc}
D. R. McMillan, Jr., {\it Some contractible open 3-manifolds}, Trans. 
Amer. Math. Soc. 102 (1962), 373--382. 

\bibitem{MSY}
W. Meeks, L. Simon, S. T. Yau, {\it Embedded minimal surfaces, exotic spheres,
and manifolds with postive Ricci curvature}, Annals of Math, 116 (1982), 621--659.

\bibitem{Me} G. Mess, {\it Centers of 3-manifold groups and groups which are 
coarse quasi-isometric to planes}, preprint. 

\bibitem{My homology}
R. Myers, {\it Homology cobordisms, link concordances, and hyperbolic
3-manifolds}, Trans. Amer. Math. Soc., 278 (1983), 271--288.


\bibitem{My genus one}
R. Myers, {\it Contractible open 3-manifolds which are not covering spaces},
Topology, 27 (1988), 27--35.

\bibitem{My excel}
R. Myers, {\it Excellent 1-manifolds in compact 3-manifolds}, Topology Appl. 49 (1993),
115--127.


\bibitem{My attach}
R. Myers, {\it Attaching boundary planes to irreducible open 3-manifolds},
Quart. J. Math. Oxford Ser. (2), to appear. 

\bibitem{My endsum}
R. Myers, {\it End sums of irreducible open 3-manifolds}, Oklahoma State 
University Mathematics Department Preprint Series (1996). 


\bibitem{My cover}
R. Myers, {\it Contractible open 3-manifolds which non-trivially cover only
non-compact 3-manifolds}, Oklahoma State University Mathematics Department 
Preprint Series (1996).  


\bibitem{Sh}
P. Shalen, {\it Dendrology of groups: an introduction}, Essays in Group 
Theory, Math. Sci. Res. Inst. Publ. 8, Springer, New York--Berlin (1987), 
265--319. 

\bibitem{St}
J. Stallings, {\it On torsion free groups with infinitely many ends}, Annals of Math. 88 (1968),
312--334.


\bibitem{Ti-Wr}
F. Tinsley and D. Wright, {\it Some contractible open manifolds and coverings of
manifolds in dimension three}, Topology Appl., to appear.

\bibitem{Wa}
F. Waldhausen, {\it On irreducible 3-manifolds which are sufficiently large},
Ann. of Math., 87 (1968), 56--88.

\bibitem{Wl}
C. T. C. Wall, {\it Open 3-manifolds which are 1-connected at infinity}, 
Quart. J. Math. Oxford Ser. (2) 16 (1965), 263--268. 

\bibitem{Wi}
B. Winters, {\it Properly homotopic nontrivial planes are parallel}, Topology Appl. 48 (1992),
235--243.

\bibitem{Wr}
D. Wright, {\it Contractible open manifolds which are not covering spaces},
Topology 31 (1992), 281--291.



\end{thebibliography}
\end{document}